\documentclass[12pt]{article}
\usepackage{amsmath,amsthm,amssymb,amsfonts}

\newtheorem{theorem}{Theorem}[section]
\newtheorem{proposition}[theorem]{Proposition}
\newtheorem{lemma}[theorem]{Lemma}
\newtheorem{claim}[theorem]{Claim}
\newtheorem{remark}[theorem]{Remark}
\newtheorem{corollary}[theorem]{Corollary}

\DeclareMathOperator{\wt}{wt}

\DeclareMathOperator{\tr}{tr}
\DeclareMathOperator{\End}{End}
\DeclareMathOperator{\Res}{Res}
\DeclareMathOperator{\id}{id}
\DeclareMathOperator{\ad}{ad}
\DeclareMathOperator{\Span}{Span}
\DeclareMathOperator{\Ker}{Ker}

\begin{document}
\title{The space of graded traces for holomorphic VOAs of small central 
charge\thanks{The majority of this work was part of the author's doctoral 
dissertation, which was completed at the University of California, Santa Cruz 
under the direction of Professor Geoffrey Mason.}}
\author{Katherine L. Hurley}
\date{\today}

\maketitle

\begin{abstract}
It is one of the remarkable results of vertex operator algebras (VOAs) that 
the graded traces (one-point correlation functions) of holomorphic VOAs 
are modular functions.  This paper explores the question 
of which modular functions arise as the graded traces of holomorphic VOAs.  
For VOAs of small central charge, i.e., $c\le 24$, and a non-zero weight-one
subspace we find that the only 
conditions imposed on the modular functions are those that arise easily out 
of our condition that the VOAs be of CFT type, that is that they have 
no negative-weight subspaces and their zero-weight subspace is generated by 
the vacuum vector.
\end{abstract}

\section{Introduction}
One of the striking features of rational vertex operator algebras (VOAs)
is the modularity of their graded traces, also called one-point correlation
functions.    
Zhu proved that for any holomorphic VOA, $V$, satisfying a finiteness 
condition, and certain $v\in V$, 
the graded trace of $v$ , $Z(v,\tau)$, is a meromorphic modular form, possibly 
with a character~\cite[Theorem 5.3.2]{Zhu:mi}. 
We will refer to this result as Zhu's 
Theorem, although Dong, Li and Mason weakened Zhu's assumptions to those 
stated~\cite{DLM:mi}~\cite{DLM:tsmm}.   
In this paper, we investigate an inverse to Zhu's Theorem.  
For a fixed VOA $V$, and modular form $f$, 
does there exist a $v\in V$ such that $Z(v,\tau)=f$?  
For the relevant definitions, please see Section~\ref{Definitions} 

We consider strongly rational holomorphic VOAs of small central 
charge $c$, i.e., $c\le 24$. 
We make use of work of Dong and Mason on the classification of such 
VOAs~\cite{DM:hVOAscc}.  
It is well known that for holomorphic VOAs, $c$ is positive and 
$8|c$.  Thus 
we are considering the cases $c=8$, $c=16$ and $c=24$.  Dong and Mason 
prove that for $c=8$ and $c=16$, the only strongly holomorphic VOAs are the
lattice VOAs generated by the $E_8$ root lattice for $c=8$ and the lattices
$E_8+E_8$ and $D_{16}^+$ for $c=16$.  

However, for the $c=24$ case Dong and Mason do not achieve a 
complete 
classification.  They do show that the weight-one subspace $V_1$, which
forms a Lie algebra, is either 
zero, abelian or semisimple.  In this paper we treat the $V_1$ abelian and 
$V_1$ semisimple cases. Note that we do not use Schellekens classification 
of $c=24$ holomorphic VOAs, which provides a list of 71 possible Lie algebras 
for $V_1$~\cite{Schell}, as his classification is not mathematically 
rigorous.  

Frenkel, Lepowsky and Meurmann  conjecture that the 
only holomorphic $c=24$ VOA with $V_1=0$ is the moonshine module~\cite{FLM}.  
If this conjecture holds, then the current paper and Dong and Mason's paper 
\emph{Monstrous moonshine of higher weight} find the space of graded traces 
for all strongly rational holomorphic VOAs of small central 
charge~\cite{DM:MM}. 

In this work we make the assumptions about our VOA $V$ used in Zhu's 
Theorem, i.e., holomorphic and $C_2$-cofinite, which are explained in 
Section 2.  Additionally, we assume that the VOAs are of CFT-type (here
CFT stands for conformal field theory).  This assumption places immediate 
restrictions on the set of modular forms  achieved as graded traces of 
the elements of $V$.  The main result of this paper is that
for small central charge these are the only restrictions.

Let $M_k$ denote the space of holomorphic modular forms of weight $k$ and let 
$\eta(\tau)$ 
be the Dedekind eta function; see equation~\eqref{defn_eta}. 
\begin{theorem}
\label{main_c=8,16}
Let $V$ be a strongly rational holomorphic vertex operator algebra with 
central charge $c$ equal to $8$ or $16$.  
For any holomorphic modular form $f(\tau)$ with $\wt f\ge c/2$,
there exists an $v\in V$ such that 
\[Z(v,\tau)=\frac{f(\tau)}{\eta(\tau)^c}.\]
\end{theorem}
This means that every possible modular form is attained as the graded trace
of some element of $v$.  Combining this theorem with Zhu's Theorem 
characterizes the space of graded traces.
\begin{corollary}
\label{corollary_c=8,16}
Let $V$ be a strongly rational holomorphic vertex operator algebra 
with central charge $c$ equal to 8 or 16.  
The space of graded traces of $V$ is 
\[\dfrac{1}{\eta(\tau)^{c}}\bigoplus_{k\ge c/2}M_k.\]
\end{corollary}

For $c=24$ we attain the following theorem.
\begin{theorem}
\label{main_c=24}
Let $V$ be a strongly rational holomorphic vertex operator algebra 
with central charge $c=24$ and $V_1\ne 0$.  
Then for any holomorphic modular form $f(\tau)$ with $\wt f \ge 14$, 
there exists a vector $v$ in $V$ such that 
\[Z(v,\tau)= \frac{f(\tau)}{\Delta(\tau)}.\] 
\end{theorem}
Here $\Delta(\tau)=\eta(\tau)^{24}$ is the discriminant function.
With Zhu's Theorem, this describes the space of graded traces for $v$ with 
square-bracket weight greater than 
zero.  It remains to consider the vacuum vector, $\boldsymbol{1}$.
In the $c=24$ case, $Z(\boldsymbol{1},\tau)$ is a modular function with 
a pole of order $1$ at infinity and no poles in the upper half plane. 
The modular functions are rational 
functions in the modular invariant $J(\tau)$.  If normalize
the $q$-expansion of $J(\tau)$ to be $q^{-1}+0+196884q+\cdots$, then 
$Z(\boldsymbol{1},\tau)=\dim V_1 + J(\tau)$, and we can describe the 
space of graded traces of $V$.
\begin{corollary}
\label{corollary_c=24}
Let $V$ be a strongly rational holomorphic vertex operator algebra 
with central charge $24$ and $V_1\ne 0$. 
The space of graded traces of V is 
\[ \mathbb{C}(\dim V_1 + J(\tau))\oplus \frac{1}{\Delta(\tau)}
\bigoplus_{k\ge 14}M_k.\]
\end{corollary}

\section{Definitions}
\label{Definitions}
Let $(V,Y,\omega,\boldsymbol{1})$ be a vertex operator algebra.  
We recall the portions of the definition of a VOA relevant to this paper. For
the full axioms, see Frenkel, Huang and Lepowsky~\cite{FHL}.  The vector 
space $V$ is $\mathbb{Z}$-graded: $V=\bigoplus_{n\ge n_0}V_n$ and elements
of $V_n$ are called homogeneous of weight $n$.  We
denote the vertex operator $Y(v,z)=\sum_{n\in\mathbb{Z}}v(n)z^{-n-1}$.  
Here each $v(n)$ is an endomorphism of $V$.
The most important VOA identity is the Jacobi identity, which we do not 
state here; however,
in this paper, we will make frequent use of associativity,
\begin{equation}
\label{associativity}
(u(m)v)(n)=\sum_{i\ge 0}(-1)^i\binom{m}{i}
\big(u(m-i)v(n+i)-(-1)^mv(m+n-i)u(i)\big),
\end{equation}
which follows directly from the Jacobi identity.

There is a distinguished element $\omega\in V_2$, called the conformal 
vector.  Its modes are denoted
by $Y(\omega,z)=\sum_{n\in\mathbb{Z}}L(n)z^{-n-2}$.  Clearly 
$L(n)=\omega(n+1)$.  The redundant notation is introduced as the $L(n)$ 
operators make $V$ a module for the Virasoro Lie algebra.  That is
\[ [L(m),L(n)]=(m-n)L(m+n)+\frac{m^3-m}{12}\delta_{m+n,0}\id c.\]
The constant $c$ is called the \emph{central charge} of $V$.  The $L(0)$ 
operator determines the grading of $V$, i.e., for $v\in V_n$, $L(0)v=nv$.  
The highest-weight vectors for this representation of the Virasoro algebra, 
i.e., those $v\in V$ such that $L(n)v=0$ for all $n>0$, 
are called the \emph{highest-weight vectors} of $V$.
   
For homogeneous $v\in V_n$, the \emph{zero mode} of $v$ is defined 
$o(v):=v(n-1)$.  
By extending linearly, we can define $o(v)$ for all $v\in V$.  The zero mode
preserves the grading of $V$.  The graded 
trace of $v$ on the VOA is defined by:
\[
Z(v,q)=q^{-c/24}\sum_{n\in\mathbb{Z}}\tr\big|_{V_n}o(v) q^n.
\]
Here $q$ is just a formal variable; however, later we will consider $q$ as 
a complex variable of modulus less than one.  As usual in the theory of 
modular forms, we  define $\tau$ in the 
complex upper half plane by $q=e^{2\pi i \tau}$.  Hence 
we also use $Z(v,\tau)$ to denote the graded trace. 

Let $V$ be a vertex operator algebra.  We say that $V$ is \emph{rational}, 
if every admissible module for $V$ is completely reducible.~\cite{DLM:trVOA}  
It is \emph{holomorphic}, if it is rational and has only one irreducible 
module.  Furthermore, we say that $V$ is \emph{$C_2$-cofinite}, if the subspace
generated by $\{u(-2)v|u,v\in V\}$ is of finite codimension in $V$.  It is 
conjectured that rational implies $C_2$-cofinite.

A VOA is of \emph{CFT-type}, if there are no negatively graded 
subspaces and the zero-graded space is one dimensional and generated by the
vacuum vector, i.e., 
\[ V=\mathbb{C}\boldsymbol{1}+V_1+V_2+\cdots.\] 
To be of \emph{strong CFT-type}, we additionally require that $L(1)V_1=0$.  
For holomorphic VOAs, CFT-type implies strong CFT-type~\cite{DM:rVOAecc}.
Following the terminology of Dong and Mason,
we call rational $C_2$-cofinite VOAs of strong CFT-type \emph{strongly
rational}.  This paper concerns strongly rational holomorphic VOAs with
central charge $c\le 24$.

To study the graded traces of a VOA $(V,Y(\ ,z),\omega,\boldsymbol{1})$, Zhu 
introduced another VOA 
$(V,Y[\ ,z],\tilde{\omega},\boldsymbol{1})$, which is isomorphic to 
$(V,Y(\ ,z),\omega,\boldsymbol{1})$~\cite{Zhu:mi}.
While the underlying vector space and the vacuum vector are the same as the 
original VOA, the vertex operator is defined $Y[v,z]:= Y(v,e^z-1)e^{z \wt v}$, 
and the conformal vector is 
$\tilde{\omega}=\omega -\frac{c}{24}\boldsymbol{1}$. This induces a different
grading $V=\sum_{n\ge n_0}V_{[n]}$. We denote the square-bracket vertex 
operators by $Y[v,z]=\sum_{n\in\mathbb{Z}}v[n]z^{-n-1}$.
Let $v\in V_m$. To expand the square-bracket operator $v[n]$
in terms of round-bracket operators and vice-versa we use:
\begin{align}
\label{sqr-rnd} 
v[n]&=\Res_z Y(v,z)(\log (1+z))^n(1+z)^{m-1},\\
\label{rnd-sqr}
v(n)&=\Res_z Y[v,z](e^z-1)^ne^{z(1-m)}    
\end{align}
which follow from the definition of $Y[v,z]$.  Using these formulas, we can
show that $\bigoplus_{n\le N}V_{[n]}=\bigoplus_{n\le N}V_n$.

\section{Modular Forms and Zhu's Theorem}
Throughout this section $V^c$ will denote a holomorphic, $C_2$-cofinite 
vertex operator algebra of CFT-type with central charge $c$.  
For $c=24$, we also assume that
the weight-one homogeneous space is non-zero, i.e., $V^{24}_1\ne 0$.  

Let $v\in V^c_{[k]}$.  
From Zhu we know that $Z(v,q)$ the graded trace 
of $v$ converges to $q^hf(q)$ for some $f$ a holomorphic function on 
$D=\{q\in\mathbb{C}\, |\, |q|<1\}$~\cite[Theorem 4.4.1]{Zhu:mi}. 
Additionally, Zhu shows that the graded trace of 
$v$ is a modular form of weight $k$ in the following sense:  
Change variable to $\tau$, via $q=e^{2\pi i \tau}$.  
then for
all $v\in V_{[k]}$, $g=\left(\begin{smallmatrix}a&b\\c&d\end{smallmatrix}
\right)\in PSL_2(\mathbb{Z})$
and $\tau\in \mathbb{H}$
\begin{equation}
\label{g_action}
Z(v,g\cdot\tau)=\psi(g)(c\tau+d)^kZ(v,\tau),
\end{equation}
where $g\cdot \tau$ is the usual action of $g$ by linear fraction 
transformations, and $\psi:PSL_2(\mathbb{Z})\rightarrow \mathbb{C}^*$ is a 
group 
character~\cite[Theorem 5.3.3]{Zhu:mi}.

For this paper, in addition to $V^c$ holomorphic and $C_2$-cofinite we
assume that $V^c$ is of CFT-type.  So
\begin{equation}
\label{CFT_grad_trace}
Z(v,q)=q^{-c/24}\sum_{n\ge 0}\tr\big|_{V_n} o(v)q^n.
\end{equation}
Thus there exists an $f$, holomorphic on $D$, such that
\begin{equation*}
Z(v,q)=q^{-c/24}f(q).
\end{equation*}
Furthermore, because the abelianization of $PSL_2(\mathbb{Z})$ is a cyclic
group of order six generated by the coset of  
$T=\left(\begin{smallmatrix}1&1\\0&1\end{smallmatrix}\right)$, 
equation~\eqref{CFT_grad_trace} allows 
us to determine the character in equation~\eqref{g_action}.  
Let $\chi$ be the character of 
$PSL_2(\mathbb{Z})$ determined by $\chi(T)=e^{2\pi i/6}$, then
\begin{equation*}
Z(v,g\cdot \tau)=\chi(g)^{-c/4}(c\tau+d)^kZ(v,\tau),  
\end{equation*}
for all $v\in V_{[k]}$.

Let $P_k^{-c/24}$ denote the space of functions $h$ on $D$ such that there 
exists a function $f$ holomorphic on $D$ with 
\begin{align}
\label{q-expan}
h(q)&=q^{-c/24}f(q)\\
\label{trans}  
h(g\cdot \tau)&=\chi(g)^{-c/4}(c\tau+d)^kh(\tau),
\end{align}
for all $g\in PSL_2(\mathbb{Z})$ and $\chi$ as above.
We refer to the elements of $P_k^{-c/24}$ as
weight-$k$ modular forms with a pole of order $c/24$ at infinity.
Let $Z(W):=\{Z(v,q)|v\in W\}$ be the space of graded traces for $W$.
Using this notation, Zhu's Theorem implies that 
$Z(V^c_{[k]})\subseteq P_k^{-c/24}$ for all integers $k$. 
   
Note that the usual definition of a modular form $f$ of weight $k$ requires 
that $f$ be holomorphic on $D$ and have trivial character.  We will call such 
forms \emph{holomorphic modular forms} of weight $k$ and denote their 
space by $M_k$.
We recall some basic results about holomorphic modular forms.  For proofs and
further discussion see Serre's \emph{A Course in Arithmetic} or another 
introductory text on modular forms~\cite{Ser:CA}.  

There are no non-zero holomorphic modular forms with odd or negative integral 
weight, and the only weight-zero forms are the constant functions.  There are
also no non-zero holomorphic modular forms of weight two, but for every even 
integer $k\ge 4$, there exists a form, $G_k$, called the Eisenstein series of 
weight-$k$.  The Eisenstein series can be defined by their expansions at 
infinity.
\begin{equation}
G_k(q):= -\frac{B_k}{k!}+\frac{2}{(k-1)!}\sum_{n\ge 1}\sigma_{k-1}(n)q^n.
\end{equation}
Here $B_k$ is the $k$th Bernoulli number defined by 
$\frac{t}{e^t-1}=\sum_{k\ge 0}B_k\frac{t^k}{k!}$, and 
$\sigma_\ell(n):=\sum_{d|n}d^\ell$.

A key subspace of $M_k$ is the space of cusp forms 
\[
M_k^0:=\Big\{\sum_{n\ge 0} a_n q^n\in M_k\Big|a_0=0\Big\}.
\]
The smallest integer $k$ such that $\dim M_k^0>0$ is 12, and $M_{12}^0$ is 
spanned by the discriminant, 
$\Delta:=10,800\left(20(G_4)^3-49(G_6)^2\right)$, which has a simple zero at
$\tau=\infty$ ($q=0$) and no other zeros or poles.  The discriminant has 
particularly nice product formula:
\begin{equation*}
\Delta(q)=q\prod_{n=1}^\infty (1-q^n)^{24},
\end{equation*} 
and multiplication by the discriminant is a linear isomorphism of $M_k$ and 
$M_{k+12}^0$.

In general multiplication of two forms of weights $k$ and $\ell$ gives 
a form of weight $k+\ell$.  This makes the space of all holomorphic modular
forms $M:=\oplus_{k\ge 0} M_k$ a graded algebra.  Indeed, it is a 
polynomial algebra generated by $G_4$ and $G_6$.  

The $24$th root of the discriminant is the  
Dedekind eta function 
\begin{equation}
\label{defn_eta}
\eta(q):=q^{1/24}\prod_{n>0}(1-q^n).
\end{equation} 
(See for example Lang's book on modular forms~\cite{Lang:IMF}.)
The elements $T=\left(\begin{smallmatrix}1&1\\0&1\end{smallmatrix}\right)$ and
$S=\left(\begin{smallmatrix}0&-1\\1&0\end{smallmatrix}\right)$ generate the 
modular group, $PSL_2(\mathbb{Z})$.
The eta function is a modular form of weight $1/2$ and multiplier system
$\phi$ defined by $\phi(T)=e^{2\pi i/24}$ and $\phi(S)=\sqrt{\frac{1}{i}}$, 
i.e.,
for all $g=\left(\begin{smallmatrix}a&b\\c&d\end{smallmatrix}\right)
\in PSL_2(\mathbb{Z})$,
\begin{equation}
\label{eta_trans} 
\eta(g\cdot\tau)=\phi(g)(c\tau+d)^{1/2}\eta(\tau). 
\end{equation}
In the same way that multiplication by the discriminant gives a linear 
isomorphism of $M_k$ and $M_{k+12}^0$, multiplication by $\eta(\tau)^c$ gives 
a linear isomorphism of $P_{k}^{-c/24}$ and $M_{k+c/2}$.

\begin{proposition}
\label{P_and_M_iso}
Let $8|c$, multiplication by $\eta(\tau)^c$ is a linear isomorphism of 
$P^{-c/24}_k$ and $M_{k+c/2}$. 
\end{proposition}
\begin{proof}
Let  $h(\tau)$ be in $P_k^{-c/24}$.  Since $h(\tau)$ and $\eta(\tau)$ are  
holomorphic for all finite $\tau$,
the product $\eta(\tau)^c f(\tau)$ is holomorphic except 
possibly at infinity. From the 
$q$-expansions, equations~\eqref{q-expan} and~\eqref{defn_eta}, 
the product is also holomorphic at infinity.  From equations~\eqref{trans} 
and~\eqref{eta_trans}, for any $g\in PSL_2(\mathbb{Z})$,  
\[\eta(g\cdot \tau)^ch(g\cdot \tau)=\phi(g)^c\chi(g)^{-c/4}(c\tau+d)^{k+c/2}
\eta(\tau)^ch(\tau).\]
We compute that $\phi(T)^c\chi(T)^{-c/4}=1$. Since $8|c$, $\phi(S)^c=1$, 
and since $\chi(S)=1$, $\chi(S)^{-c/4}=1$.  So $\phi(g)^c\chi(g)^{-c/4}=1$ for
all $g\in PSL_2(\mathbb{Z})$
Thus $\eta(\tau)^ch(\tau)$, is in $M_{k+c/2}$.

It is clear that multiplication by $\eta(\tau)^c$ is linear.  Moreover, because
$\eta(\tau)$ has 
$q$-expansion $\eta(q)=q^{1/24}f(q)$, such that $f(q)$ has no zero or 
poles in $D$, division by $\eta(\tau)^c$ maps $M_{k+c/2}$ back to 
$P_k^{-c/24}$.  
Thus multiplication by $\eta(\tau)^c$ is invertible, and an isomorphism.
\end{proof}

Using this proposition, Zhu's Theorem implies that 
\[Z(V_{[k]}^c)\subseteq \frac{1}{\eta(\tau)^c}M_{k+c/2}.\]
Since assuming that $V^c$ is holomorphic implies the central charge $c$ is an 
integer divisible by $8$, $k+c/2$ is an integer.   
Because $V^c$ is of CFT-type, $Z(V^c)\subseteq \bigoplus_{k\ge 0}P_k^{-c/24}$. 
Thus 
\begin{equation}
\label{Zhu_restated}
Z(V^c)\subseteq \frac{1}{\eta(\tau)^c}\bigoplus_{k\ge c/2}M_k.
\end{equation} 

For $c=8$ and $c=16$ our main theorem, Theorem~\ref{main_c=8,16}, says that all
such modular forms are achieved as the graded traces of elements of $V^c$.
So our theorem and Zhu's Theorem imply Corollary~\ref{corollary_c=8,16},
i.e.,    that
\[Z(V^c)= \frac{1}{\eta(\tau)^c}\bigoplus_{k\ge c/2}M_k.\]

For $c=24$ and $k=0$, Zhu's Theorem says that 
$Z(V_{[0]}^{24})\subseteq \frac{1}{\Delta(\tau)}M_{12}$.  In this case,
CFT-type imposes another another restriction
on the image of $V_{[0]}^{24}$.
Since $\dim V^{24}_{[0]}=1$ and $\dim M_{12}=2$, 
$Z( V^{24}_{[0]})$ is a proper subspace of $\frac{1}{\Delta(\tau)}M_{12}$. 
It is generated by the graded trace of the vacuum vector, $\boldsymbol{1}$.
Thus 
\begin{equation}
\label{Zhu's_Theorem_c=24}
Z(V^{24})\subseteq \mathbb{C}Z(\boldsymbol{1},\tau)
\oplus \frac{1}{\Delta(\tau)}\bigoplus_{k\ge 14} M_k.
\end{equation}
Our main theorem for $c=24$, Theorem~\ref{main_c=24}, says that
$\frac{1}{\Delta(\tau)}\bigoplus_{k\ge 14}M_k\subseteq Z(V^{24})$. So 
equation~\eqref{Zhu's_Theorem_c=24} is  actually an equality.  This is 
essentially Corollary~\ref{corollary_c=24}.  For the corollary, note that 
from Proposition~\ref{P_and_M_iso},
$Z(\boldsymbol{1},\tau)\in P^{-1}_{[0]}=
\frac{1}{\Delta(\tau)}M_{12}$, which has basis the constant function $1$ and
the modular invariant $J(q)=q^{-1}+196884q+\cdots$.  Moreover 
$Z(\boldsymbol{1},q)=q^{-1}+\dim V_1 + \dim V_2 q +\cdots$.  So 
$Z(\boldsymbol{1},\tau)=\dim V_1 + J(\tau)$, and we get
Corollary~\ref{corollary_c=24} exactly.

\section{Proof of Main Theorems}
To prove our main theorems, we will use several results that Dong and Mason 
developed to determine 
the space of graded traces of the moonshine module. 
They state their results for the moonshine module, but the same proofs work 
for the more general statements we give here.
\begin{proposition}[\cite{DM:MM}]
\label{L[-2]_traces}
Let $V$ be a strongly rational holomorphic VOA.  For any $\ell\ge 0$,
the graded trace
of $L[-2]^\ell\boldsymbol{1}\in V_{[2\ell]}$ is non-zero.  More specifically,
\[Z(L[-2]^\ell\boldsymbol{1},\tau)=q^{-c/24}\sum_{n\ge 0}a_nq^n,\]
such that $a_0\ne 0$.
\end{proposition}

The above proposition implies that 
$\eta(\tau)^c Z(L[-2]^\ell\boldsymbol{1},\tau)$ 
is in $M_{2\ell+c/2}$ but not in $M^0_{2\ell+c/2}$.  
Since $\dim M_k/M_k^0\le 1$, this means that 
for all integers $k\ge c/2$,
\begin{equation}
\label{M_k_decomp}
M_k=\mathbb{C}\eta(\tau)^c Z(L[-2]^{(2k-c)/4}\boldsymbol{1},\tau)\oplus M_k^0.
\end{equation}
We have restricted our problem to 
determining for which cusp forms $f(\tau)\in M_k^0$, $k\ge c/2$, there exists 
a $v\in V^c$ such that $\eta(\tau)^c Z(v,\tau)=f(\tau)$.
 
For $c=8$, this question has already been answered. 
The unique strongly rational holomorphic VOA with central 
charge $8$ is the $E_8$ 
root lattice VOA, $V_{E_8}$~\cite{DM:hVOAscc}.  For this VOA, Dong Mason and
Nagatomo use a result of Waldspurger to show that for all cusp forms $f(\tau)$,
there exists a highest-weight vector $v$ such that 
$\eta(\tau)^8 Z(v,\tau)=f(\tau)$~\cite{DMN}.
Thus $m_8\circ Z$ maps $V_{E_8}$ onto $\bigoplus_{k\ge 4} M_k$.  This proves
Theorem~\ref{main_c=8,16} for $c=8$.

For the remaining central charges we use the following key lemma.
\begin{lemma}[main lemma]
\label{traces_lemma}
Let $V$ be a strongly rational holomorphic VOA with central charge 16 or 24 
and $V_1\ne 0$. Then there exists a $v\in V_{[4]}$ such that 
\begin{align*}
\tr o(v)\big|_{V_0}&=0,\\
\tr o(v)\big|_{V_1}&\ne 0.
\end{align*}
\end{lemma}

The proof of Lemma~\ref{traces_lemma} is quite lengthy and it is postponed to 
Section~\ref{proof_of_traces_lemma}. The lemma is used to prove the following 
proposition.
\begin{proposition}
\label{smallest_cusps}
Let $V^c$ be a strongly rational holomorphic VOA with central charge
$c$ equal to 16 or 24 and $(V^c)_1\ne 0$.  
There exists a $v\in (V^{16})_{[4]}$ such that 
\[ Z(v,\tau)=\eta(\tau)^{8},\]
and there exists a $v\in (V^{24})_{[4]}$ such that
\[ Z(v,\tau)=G_4(\tau).\]
\end{proposition} 
\begin{proof}
Assuming Lemma~\ref{traces_lemma} gives the existence of a $v\in V^c$ such that
$\eta(\tau)^cZ(v,\tau)$ is a non-zero element of $M^0_{4+c/2}$. 
For $c=16$, this is $M^0_{12}$, which is a one-dimensional space generated
by $\Delta(\tau)=\eta(\tau)^{24}$.    
Thus, adjusting by a constant multiple if necessary, 
$\eta(\tau)^{16} Z(v,\tau)=\Delta(\tau)$, i.e.
\[ Z(v,\tau)=\frac{\Delta(\tau)}{\eta(\tau)^{16}}=\eta(\tau)^8.\]
 
For $c=24$, we have that $\Delta(\tau)Z(v,\tau)$ is a non-zero element 
of $M_{16}^0$, which is also one dimensional and generated by 
$\Delta(\tau)G_4(\tau)$.  So adjusting by a constant if necessary, 
\[Z(v,\tau)=\frac{\Delta(\tau)G_4(\tau)}{\eta(\tau)^{24}}=G_4(\tau).\]
\end{proof}

We do not need prove a lemma like Lemma~\ref{traces_lemma} for $v\in V_{[k]}$,
$k>4$. For those vectors we use a result of Dong and Mason which 
describes the space of graded traces for the Virasoro module generated by a 
highest-weight vector $w$.

To describe this space of graded traces, we use a derivation
$\delta_k:P_k^{-c/24}\rightarrow P_{k+2}^{-c/24}$.  This derivation is a 
generalization of the derivation $\delta_k:M_k\rightarrow M_{k+2}$ described
in Lang.~\cite{Lang:IMF}[Chapter X, Section 5]

Define the ``Eisenstein series'' $G_2(\tau)$ by 
\[G_2(q)=-\frac{1}{12}+2\sum_{n=1}^\infty\sigma_1(n)\, q^n.\]
For $g=\left(\begin{smallmatrix}a&b\\c&d\end{smallmatrix}\right)
\in PSL_2(\mathbb{Z})$,
\[G_2(g\cdot\tau)=(c\tau+d)^2 G_2(\tau)-\frac{c}{2\pi i}(c\tau+d).\]
Thus $G_2$ is not a modular form.  For $f(q)\in P_k^{-c/24}$, define 
\[\delta_k f(q) =q\frac{d}{dq}f(q)+k G_2(q) f(q).\]
\begin{proposition}
The map $\delta_k$ takes $P_k^{-c/24}$ to $P_{k+2}^{-c/24}$ and is 
a derivation in the sense that if $g\in M_k$ and 
$f\in P^{-c/24}_\ell$, then
\[\delta_{k+\ell}(gf)=\big(\delta_k(g)\big) f+
                                  g\big(\delta_\ell(f)\big).\]   
\end{proposition}
\begin{proof}
The proof relies on  $\delta_k:M_k\rightarrow M_{k+2}$ 
being a derivation of the space of holomorphic modular forms, and 
equation~\eqref{delta_on_P}, below.

The logarithmic derivative of $\eta(\tau)$ is $-\pi i G_2(\tau)$ 
and $q\frac{d}{dq}=\frac{1}{2\pi i}\frac{d}{d\tau}$.   
Using these facts, the definition of $\delta_k$ and the quotient rule, 
we show that 
if $f(\tau)\in P_k^{-c/24}$ equals $\frac{h(\tau)}{\eta(\tau)^c}$ with
$h(\tau)\in M_{k+c/2}$, then
\begin{equation}
\label{delta_on_P}
\delta_k (f(\tau))=\frac{\delta_{k+c/2}(h(\tau))}{\eta(\tau)^c}.
\end{equation}
Since $\delta_{k+c/2}(h(\tau))\in M_{k+c/2+2}$, 
$\delta_k(f(\tau))\in P^{-c/24}_{k+2}$.
To see that $\delta$ has the derivation property given in the 
proposition, simply
use equation~\eqref{delta_on_P} and that 
$\delta$ is a graded derivation on the space of holomorphic modular forms.
 \end{proof}  
Define a $\delta$-submodule of 
$P^{-c/24}$ to be a linear subspace closed under multiplication by elements of
$M$ and the action of $\delta_k$. 
\begin{theorem}[\cite{DM:MM}]
\label{delta_module}
Let $V$ be a $C_2$-cofinite holomorphic VOA and let $w$ be a highest-weight
vector of positive weight.  The space of graded traces consisting of all 
$Z(v,\tau)$ for $v$ in the Virasoro module 
generated by $w$ is the $\delta$-module generated by $Z(w,\tau)$.
\end{theorem}

Thus if there exists a highest-weight vector 
$w\in V^{24}$ with graded trace equal to $G_4(\tau)$,  
then $Z(V^{24})$ contains the $\delta$-module 
generated by $G_4(\tau)$.  The 
$\delta$-module generated by $G_4(\tau)$ also contains $G_6(\tau)$, 
as $\delta_4 G_4= 14 G_6$~\cite{Lang:IMF}.  
Since $G_4$ and $G_6$ generate the space of 
holomorphic forms as a polynomial algebra, the $\delta$-module generated by 
$G_4(\tau)$ contains all the holomorphic modular forms except 
the constant functions, which implies that 
\[
\bigoplus_{k>0}M_k \subseteq Z(V^{24}).
\]  
In light of equation~\eqref{M_k_decomp}, this would complete the proof of 
the main theorem for the $c=24$ case, Theorem~\ref{main_c=24}.

If there exists a highest-weight vector $w\in V^{16}$ with graded trace equal 
to $\eta(\tau)^8$, then $Z(V^{16})$ contains the $\delta$-module in $P^{-2/3}$
generated by $\eta(\tau)^8$.  Thus $m_{16}\circ Z(V^{16})$ contains the 
image of that $\delta$-module under $m_{16}$, which 
contains $\eta(\tau)^{16}\eta(\tau)^8=\Delta(\tau)$.  
The image of a $\delta$-module in $P^{-2/3}$ under $m_{16}$ is a 
$\delta$-module in $M$.   
The $M$-module generated by $\Delta$ in $M$ is $M^0$, 
the space of cusp forms.  So
$M^0\subseteq m_{16}\circ Z(V^{16})$.  Combining this with 
equation~\eqref{M_k_decomp} proves Theorem~\ref{main_c=8,16} 
in the $c=16$ case.

Note that while Proposition~\ref{smallest_cusps} gives the existence of 
vectors $v_1\in V^{24}$ and $v_2 \in V^{16}$ with $Z(v_1,\tau)=G_4(\tau)$ and 
$Z(v_2,\tau)=\eta(\tau)^8$, the above argument requires the existence of 
\emph{highest-weight} vectors with those graded traces.  To fill this hole in
our argument we need the following lemma.

\begin{lemma} 
\label{highest-weight_lemma}
Let $V^c$ be a strongly rational holomorphic VOA with central charge $c = 16$ 
or $24$.  If there exists a $v\in V_{[4]}$ such that $\tr o(v)|_{V_0}=0$, 
then there exists a highest-weight vector $w$ in $V^c$ such that 
$Z(w,\tau)=Z(v,\tau)$.
\end{lemma} 

We prove this lemma in the following section and the main lemma, 
Lemma~\ref{traces_lemma},
in Section~\ref{proof_of_traces_lemma}. 
This will complete the proof of our main theorem.

\section{Proof of Lemma~\ref{highest-weight_lemma}}

By a highest-weight vector in Lemma~\ref{highest-weight_lemma}, 
we mean a round-bracket highest-weight vector,
i.e., $L(n)v=0$ for all $n>0$.  However, using equations~\eqref{sqr-rnd} 
and~\eqref{rnd-sqr} 
it is not hard to show that the set of round-bracket
highest-weight vectors in $V^c_4$ equals the set of square-bracket 
highest-weight vectors in $V^c_{[4]}$. 

To analyze $V^c_{[4]}$ we  make use of the symmetric 
invariant bilinear form on square-bracket  $V^c$.  
A bilinear form $(\ ,\ )$ on $V$ is 
said to be \emph{invariant} if 
\begin{equation}
\label{VOA_invariance}
\big(Y[u,z]v,w\big)=\big(v,Y(e^{zL[1]}(-z^{-2})^{L[0]}u,z^{-1})w\big),
\end{equation}
for all $u, v, w\in V^c$~\cite{FHL}.  
Li proved that the space of invariant forms on $V$ is isomorphic to 
$\big(\frac{(V^c)_{[0]}}{L[1](V^c)_{[1]}}\big)^*$ \cite{Li:bf}.  
Since $V^c$ is of strong CFT-type,
$\big(\frac{(V^c)_{[0]}}{L[1](V^c)_{[1]}}\big)^*$ is one-dimensional, 
and the invariant bilinear form on $V^c$ is unique up to a scalar.
Henceforth, we fix $(\ ,\ )$ to be the unique invariant form on $V^c$ so that 
$(\boldsymbol{1},\boldsymbol{1})=1$.  Because $V^c$ is holomorphic, it is 
simple, and $(\ ,\ )$ is non-degenerate.
     
Harada and Lam showed that 
$V^c_{[4]} = \Ker L[1]\big|_{V^c_{[4]}} \oplus L[-1]V^c_{[3]}$ and
$V^c_{[4]} = \Ker L[2]\big|_{V^c_{[4]}} \oplus L[-2]V^c_{[2]}$, 
both sums orthogonal with respect to the square-bracket invariant 
bilinear form on $V^c$~\cite{HL:IVo}.
Because the form is non-degenerate, it follows that
\begin{equation}
\label{V_[4]_decomp}
V^c_{[4]}= (\Ker L[1] \cap \Ker L[2]) \oplus (L[-1]V^c_{[3]}+
L[-2]V^c_{[2]}).
\end{equation}
We will show that for any 
$v\in V^c_{[4]}$ such that $\tr |_{V^c_0}o(v)=0$,  
then there  
exists an element of $\Ker L[1] \cap \Ker L[2]$ with the same graded trace.
Because $\Ker L[1] \cap \Ker L[2]$  is precisely the space of square-bracket 
highest-weight vectors in $V^c_{[4]}$, 
this will prove Lemma~\ref{highest-weight_lemma}. 

From equation~\eqref{V_[4]_decomp}, 
there exist $w\in\Ker L[1] \cap \Ker L[2]$ and 
$u\in L[-1]V^c_{[3]}+L[-2]V^c_{[2]}$ such that $v=w+u$.
Assume that $\tr |_{V^c_0} o(v)=0$.  Since $w$ is a highest-weight vector, 
$w\in V^c_4$, and from the creation axiom: 
$\tr|_{V^c_0} o(w)=0$.
Therefore 
\[
\tr\big|_{V^c_0} o(u)=0.
\]
Now let $u=L[-1]a+L[-2]b$ with $a\in V^c_{[3]}$ and $b\in V^c_{[2]}$.  
Because $Z(L[-1]a,\tau)=0$ for all $a\in V^c$ (\cite{Zhu:mi} and 
\cite{DLM:mi}), $\tr|_{V^c_0} o(L[-1]a)=0$.  Thus  
\[\tr\big|_{V^c_0}o(L[-2]b)=0.\]
Furthermore $Z(L[-1]a,\tau)=0$ implies that 
\[
Z(v,\tau)=Z(w,\tau)+Z(L[-2]b,\tau),
\]
Hence the following claim completes the proof of 
Lemma~\ref{highest-weight_lemma}.

\begin{claim} Suppose that $V^c$ is a strongly rational holomorphic 
VOA
with central charge $c=16$ or $24$.
Let $b\in V^c_{[2]}$. 
If $\tr|_{V^c_0} o(L[-2]b)=0$, then $Z(L[-2]b,\tau)=0$.
\end{claim}
\begin{proof}
Write $b$ in terms of round brackets, $b= b_2+b_1+b_0$, $b_i\in V_i$.
In the square bracket VOA, the conformal element is 
$\tilde{\omega}=\omega-\frac{c}{24}\boldsymbol{1}$.  
Using equation~\eqref{sqr-rnd}, we compute
\[
L[-2]=-\frac{c}{24}\id+L(-2)+\frac{3}{2}L(-1)+\frac{5}{12}L(0)
-\frac{1}{24}L(1)+\frac{11}{720}L(2)+\cdots.
\]
From the creation axiom, 
\[ \tr\big|_{V^c_0} o(L[-2]b)= \tr \big|_{V^c_0}(L[-2]b)_0(-1),\]
where $(L[-2]b)_0$ is the projection of $L[-2]b$ into $V^c_0$.  
Since $L(0)b_0=0$ and $L(1)b_1=0$ as $V$ is strong CFT-type, 
$(L[-2]b)_0 = -\frac{c}{24}b_0+\frac{11}{720}L(2)b_2.$
Thus 
\[\tr\big|_{V^c_0} o\Bigl(-\frac{c}{24}b_0+\frac{11}{720}L(2)b_2\Bigr)=0,\]
and because $V$ is of CFT-type, 
$-\frac{c}{24}b_0+\frac{11}{720}L(2)b_2=0$.   
Furthermore $L[2]=L(2)-(1/2)L(3)+\cdots$, so $L[2]b=L(2)b_2$.
Thus 
\begin{equation}
\label{L[2]b_equation}
L[2]b=\frac{30c}{11}b_0.
\end{equation}

We now turn to the following equation from Zhu~\cite{Zhu:mi}, which holds
for $C_2$-cofinite VOAs -- see also 
\cite[eq. 5.8]{DLM:mi}.  For $b\in V_{[2]},$
\begin{equation}
\label{L[-2]_equation}
Z(L[-2]b,\tau)=\delta Z(b,\tau) + G_4(\tau)Z(L[2]b,\tau).
\end{equation}
Here $\delta$ denotes $\delta_k$ extended linearly to all $P^{-c/24}_k$. 
Note that if the coefficient of $q^{-c/24}$ in $Z(b,\tau)$ is zero, then 
the coefficient of $q^{-c/24}$ in $\delta Z(b,\tau)$ is also zero.  
Let $\ell$ be defined 
so that $b_0=\ell\boldsymbol{1}$, then, from equations~\eqref{L[2]b_equation}
and \eqref{L[-2]_equation},
\begin{align*}
Z(L[-2]b,\tau) & = \delta Z(b_2+b_1,\tau)+\delta Z(b_0,\tau) 
                   + \frac{30c}{11}G_4(\tau) Z(b_0,\tau),\\
               & = 0 q^{-c/24}+\cdots + -\frac{\ell c}{24} q^{-c/24}+\cdots\\
               &\quad +\frac{30c}{11} \Big(\frac{1}{720}+\frac{1}{3}
                 \sum_{n>0}\sigma_3(n)q^n\Big)(\ell q^{-c/24}+\cdots).
\end{align*}
Equating the coefficients of $q^{-c/24}$ gives
\[
\tr\big|_{V_0} o(L[-2]b) = -\frac{10\ell c}{11}.
\]
Hence $\frac{10 \ell c}{11}=0$ and since $c\ne 0$, $\ell=0$, 
that is $b_0=0$, and from
equation~\eqref{L[2]b_equation}, 
$L[2]b=0$.  
Additionally, $b_0=0$ implies that 
$\tr|_{V_0} o(b)=0$. Hence $\eta(\tau)^c Z(b,\tau)$ is a cusp form of weight
$4+c/2$.  
Since $c=16$ or $24$, $\eta(\tau)^cZ(b,\tau)$ is in $M^0_{10}$ or $M^0_{14}$. 
There are no non-zero cusp forms of weight 10 or 14, so $Z(b,\tau)=0$.  Using
this and $L[2]b=0$, in equation~\eqref{L[-2]_equation}, yields 
$Z(L[-2]b,\tau)=0$ as claimed.
\end{proof}
\begin{remark}The same result and proof holds for $c=8$ as well.  However, in 
that case, $Z(v,\tau)=0$ for all $v\in V^8_{[4]}$ such that 
$\tr|_{V^8_0} o(v)=0$, so the result is inconsequential.
\end{remark}

\section{Proof of Lemma~\ref{traces_lemma}}
\label{proof_of_traces_lemma}

Let $V$ be a VOA of CFT type, then $V_1$ is a Lie algebra with an invariant 
symmetric bilinear form.  For $a$ and $b$ elements of $V_1$, the Lie bracket 
is given by $[a,b]=a(0)b$, and the 
bilinear form is given by $\langle a|b\rangle\mathbf{1}=a(1)b$, 
where $\mathbf{1}$ is the vacuum vector.  We refer to $\langle\ |\ \rangle$ 
as the Li-Zamilodichov metric. The proof that $V_1$ is a
Lie algebra depends on skew-symmetry (see~\cite{FHL}) and associativity, 
equation~\eqref{associativity}.

The operators $a(m), b(n)\in \End(V)$ make $V$ a module for  
an affine Lie algebra, i.e, they satisfy the commutation relations
\begin{equation}
\label{affine}
 [a(m),b(n)]=[a,b](m+n)+m\langle a|b\rangle \id \delta_{m+n,0}.
\end{equation}
This is proved using the Jacobi identity.

To prove the main lemma, Lemma~\ref{traces_lemma}, we will explicitly compute 
the traces of certain elements
of $V_{[4]}$ on $V_0$ and $V_1$.  For $V$ of strong CFT-type, $V_1=V_{[1]}$,
and, if $a,b,c$ and $d$ are in $V_1$, then $a[-1]b[-1]c[-1]d$ is in $V_{[4]}$.

\begin{lemma}
\label{general_traces}
Let $V$ be a strongly rational vertex operator algebra. 
Let $a$, $b$, $c$ and $d$ be elements of 
$V_1$ and $n$ be the dimension of $V_1$.
\begin{eqnarray}
\label{general_trace_V_0}
\lefteqn{\tr\Big|_{V_0} o(a[-1]b[-1]c[-1]d) =}\\
\nonumber
& &\hspace{0cm} -\frac{1}{720}\big(4\langle[a,b]|[c,d]\rangle + 5\langle [a,d]
      |[b,c]\rangle\big) +\frac{1}{1152}\sum_{\pi\in S_4}\langle \pi(a)|\pi(b)
      \rangle\langle \pi(c)|\pi(d) \rangle\\
\label{general_trace_V_1}
\lefteqn{\tr\Big|_{V_1} o(a[-1]b[-1]c[-1]d) =}\\
\nonumber
& &\hspace{0cm}-\frac{n+240}{720}\big(4\langle[a,b]|[c,d]\rangle + 5\langle 
      [a,d]|[b,c]
      \rangle\big)+\frac{n-48}{1152}\sum_{\pi\in S_4}\langle \pi(a)|\pi(b)
      \rangle\langle \pi(c)|\pi(d) \rangle\\
\nonumber
& &\hspace{.25cm}-\frac{1}{48}\sum_{\pi\in S_4}\langle \pi(a)|\pi(b)\rangle
       \kappa(\pi(c)
       ,\pi(d))+\frac{1}{6}\tr\Big|_{V_1}\sum_{\pi\in S_3}a(0)\pi\big(b(0)\big)
        \pi \big(c(0)\big)\pi \big(d(0)\big)
\end{eqnarray}
\end{lemma}
\begin{proof}
Because the graded trace is defined in terms of round-bracket homogeneous 
vectors, we need to expand $a[-1]b[-1]c[-1]d$ in terms of round brackets.  
We apply equation~\eqref{sqr-rnd},
\begin{equation*}
v[m]=\Res_z Y(v,z)(\log(1+z))^m(1+z)^{\wt v -1},
\end{equation*}
repeatedly to get
\begin{equation}
\label{square_bracket_expansion}
\begin{array}{rcl}
\lefteqn{\textstyle{a[-1]b[-1]c[-1]d=a(-1)b(-1)c(-1)d}}\hspace{.6cm}&&\\
&& +\frac{1}{2}\Big(a(0)b(-1)c(-1)d+a(-1)b(0)c(-1)d+a(-1)b(-1)[c,d]\Big)\\
&& -\frac{1}{12}\Big(a(1)b(-1)c(-1)d+a(-1)b(1)c(-1)d
   +\langle c|d\rangle a(-1)b\Big)\\
&& +\frac{1}{4}\Big(a(0)b(0)c(-1)d+a(0)b(-1)[c,d]+a(-1)[b,[c,d]]\Big)\\
&& +\frac{1}{24}\Big(a(2)b(-1)c(-1)d+a(-1)b(2)c(-1)d-a(1)b(0)c(-1)d\\
&& -a(0)b(1)c(-1)d-a(1)b(-1)[c,d]-\langle b|[c,d]\rangle a
   -\langle c|d\rangle[a,b]\\
&& +3[a,[b,[c,d]]]\Big)-\frac{19}{720}a(3)b(-1)c(-1)d\\
&& +\frac{1}{48}\Big(a(2)b(0)c(-1)d +a(2)b(-1)[c,d]
   -\langle [a,b]|[c,d]\rangle\boldsymbol{1}\Big)\\
&& +\frac{1}{144}\Big(a(1)b(1)c(-1)d+\langle c|d\rangle\langle a|b \rangle
   \boldsymbol{1}\Big).
\end{array}
\end{equation}
Thus the square-bracket weight-four vector, $a[-1]b[-1]c[-1]d$ is the sum of 
round-bracket vectors of weights zero through four.  To compute its traces on 
$V_0$ and $V_1$, we compute the traces of each round-bracket 
homogeneous component separately and add the result.  Properly, we are 
actually computing the traces of the zero mode of each vector; however, we 
usually neglect to mention ``zero mode'' explicitly. 
We start with the round-bracket 
weight-zero terms, namely: 
\begin{multline}
\label{wt_0_terms}
-\frac{19}{720}a(3)b(-1)c(-1)d+\frac{1}{48}\big(a(2)b(0)c(-1)d
+a(2)b(-1)[c,d]-\langle [a,b]|[c,d]\rangle\boldsymbol{1}\big)\\
+\frac{1}{144}\big(a(1)b(1)c(-1)d+\langle c|d\rangle\langle a|b \rangle
   \boldsymbol{1}\big).
\end{multline}

Using the affine Lie algebra commutation
relations and the invariance of the bilinear form, we can show that the 
weight-zero terms,~\eqref{wt_0_terms}, equal
\[\big(-\frac{1}{180}\langle [a,b]|[c,d] \rangle-\frac{1}{144}\langle [a,d]|
  [b,c]\rangle +\frac{1}{144}(\langle a|b \rangle\langle c|d \rangle +
\langle a|c \rangle \langle b|d\rangle +\langle a|d \rangle\langle b|c\rangle)
\big)\boldsymbol{1}.\]
The zero mode of the vacuum vector, $o(\boldsymbol{1})$, is the identity. 
So the weight-zero terms have trace on any $V_k$ equal to the constant 
in the above expression times the dimension of $V_k$. Rewriting the constant
to emphasize its symmetry, the weight-zero terms contribute 
\begin{equation}
\label{tr_0-terms}
-\frac{n}{720}\big(4\langle [a,b]|[c,d] \rangle+5\langle [a,d]|[b,c]
\rangle\big) + \frac{n}{1152}\sum_{\pi\in S_4}\langle\pi(a)|\pi(c)\rangle
\langle\pi(b)|\pi(d)\rangle
\end{equation}
to the trace of $o(a[-1]b[-1]c[-1]d)$ on $V_k$, where $n=\dim V_k$.  
In particular when $n=\dim V_1$, this is the contribution of the weight-zero
terms to $\tr\big|_{V_1}o(a[-1]b[-1]c[-1]d)$. 

From the creation axiom, if $v$ is of weight $k>0$ then 
$o(v)\boldsymbol{1}=v(k-1)\boldsymbol{1}=0$.
As $V_0=\mathbb{C}\boldsymbol{1}$, only the weight-zero terms
contribute to the trace of $o(a[-1]b[-1]c[-1]d)$ on $V_0$. 
Since $\dim V_0=1$, the first part of the lemma, 
equation~\eqref{general_trace_V_0}, follows immediately 
from equation~\eqref{tr_0-terms}. 

For a weight-one vector, $u$, $o(u)=u(0)=\ad u$, the adjoint operator.  
Because $V_1$ is reductive~\cite{DM:rVOAecc}, 
$\tr\big|_{V_1}\ad u=0$ for all $u\in V_1$.  Thus 
the weight-one terms have no contribution to the trace on $V_1$.  

We turn to the weight-two terms:
\begin{multline}
\label{2_terms}
-\frac{1}{12}\Big(a(1)b(-1)c(-1)d+a(-1)b(1)c(-1)d+\langle c|d\rangle a(-1)b
\Big)\\
+\frac{1}{4}\Big(a(0)b(0)c(-1)d+a(0)b(-1)[c,d]+a(-1)\big[b,[c,d]\big]\Big).
\end{multline}
Using the affine Lie algebra relations, we compute that
\begin{align*}
a(-1)b(1)c(-1)d &= a(-1)[[b,c],d]+\langle b|c \rangle a(-1)d +
                    \langle b|d \rangle a(-1)c,
\end{align*}
\begin{align*}
a(1)b(-1)c(-1)d &= [[a,b],c](-1)d+c(-1)[[a,b],d]+\langle a|b\rangle c(-1)d\\
                &\quad  +b(-1)[[a,c],d]+\langle a|c \rangle b(-1)d+\langle a|d
                    \rangle b(-1)c.
\end{align*} 
Each of the terms in these expansions has the form $u(-1)v$ for some 
$u$ and $v$ in $V_1$, as do the terms $\langle c|d \rangle a(-1)b$ and 
$a(-1)\big[b,[c,d]\big]$ from expression~\eqref{2_terms}.  
So their zero modes have the form $\big(u(-1)v\big)(1)$.  
We use associativity, equation~\eqref{associativity}, to compute that
\[ \big(u(-1)v\big)(1)= u(-1)v(1)+v(0)u(0)+v(-1)u(1)\]
as an operator on $V_1$.  Thus 
\[ \tr\Big|_{V_1}\big(u(-1)v\big)(1)=\kappa(u,v)+2\langle u|v\rangle.\]
For ease of notation, let $(u\! :\! v)=\kappa(u,v)+2\langle u|v\rangle$.  Note 
that $(u\! :\! v)$ is an invariant symmetric bilinear form, because both 
$\kappa(u,v)$ and $\langle u|v\rangle$ are.

There are two terms remaining in expression~\eqref{2_terms}.  One of these,
$a(0)b(-1)[c,d]$ has the form $u(0)v(-1)w$, for some $u,v,w\in V_1$. 
We expand the other:  
\[ a(0)b(0)c(-1)d= a(0)[b,c](-1)d+a(0)c(-1)[b,d],\]
and see that its terms have the same form.  The zero mode of this form 
is $\big(u(0)v(-1)w)(1)$.
\begin{align*}
\tr\Big|_{V_1}\big(u(0)v(-1)w\big)(1)&= \tr\Big|_{V_1}\big([u,v](-1)w\big)(1)+ 
\tr\Big|_{V_1}\big(v(-1)[u,w]\big)(1)\\
&= ([u,v]\! :\! w)+(v\! :\! [u,w])=0
\end{align*}
So $a(0)b(-1)[c,d]$ and $a(0)b(0)c(-1)d$ both have trace zero on $V_1$.

We now sum the traces of all the terms from expression~\eqref{2_terms} and 
get: 
\begin{multline}
\label{intermediate_2-terms}
-\frac{1}{12}\Big(-([a,d]\! :\! [b,c])-([b,d]\! :\! [a,c])+([a,b]\! :\! [c,d])
-([c,d]\! :\! [a,b])
+\langle c|d\rangle (a\! :\! b)\\+\langle b|c\rangle (a\! :\! d)
+\langle b|d\rangle (a\! :\! c)
+\langle a|c\rangle (b\! :\! d)+\langle a|d\rangle (b\! :\! c)\Big)
+\frac{1}{4}([a,b]\! :\! [c,d]).
\end{multline}
We wish to simplify this.  
Using the Lie-algebra Jacobi identity, it is easy to show that 
\begin{equation}
\label{[a,c]id}
([a,c]\! :\! [b,d])=([a,b]\! :\! [c,d])+([a,d]\! :\! [b,c]),
\end{equation}
for any invariant symmetric bilinear form.  
So expression~\eqref{intermediate_2-terms} equals
\begin{multline*}
\frac{1}{3}([a,b]\! :\! [c,d])+\frac{1}{6}([a,d]\! :\! [b,c])
-\frac{1}{12}\Big(\langle a|b\rangle(c\! :\! d)+\langle a|c \rangle(b\! :\! d)
\\
+\langle a|d\rangle(b\! :\! c)+\langle b|c \rangle (a\! :\! d)
+\langle b|d\rangle (a\! :\! c)+\langle c|d\rangle(a\! :\! b)\Big).
\end{multline*}
Substituting $\kappa(u,v)+2\langle u|v\rangle$ back in for $(u\! :\! v)$ and 
rewriting to emphasize the symmetry gives
\begin{multline}
\label{tr_2-terms}
 \frac{2}{3}\langle [a,b]|[c,d]\rangle +\frac{1}{3}\langle [a,d]|[b,c]\rangle
+\frac{1}{3}\kappa([a,b],[c,d])+\frac{1}{6}\kappa([a,d],[b,c])\\
-\frac{1}{24}\sum_{\pi\in S_4}\langle \pi(a)|\pi(b)\rangle \langle
\pi(c)|\pi(d)\rangle
-\frac{1}{48}\sum_{\pi\in S_4}\langle \pi(a)|\pi(b)\rangle \kappa(\pi(c),
\pi(d)).
\end{multline}
This is the contribution of the weight-two terms to 
$\tr\big|_{V_1}o(a[-1]b[-1]c[-1]d)$. 
  
Now let's do the weight-three terms.
These are:
\[
\frac{1}{2}\Big(a(0)b(-1)c(-1)d+a(-1)b(0)c(-1)d+a(-1)b(-1)[c,d]\Big).
\]
Expanding the first two terms this becomes
\begin{multline}
\label{3-terms_expan}
\frac{1}{2}\Big([a,b](-1)c(-1)d+b(-1)[a,c](-1)d+b(-1)c(-1)[a,d]\\
+a(-1)[b,c](-1)d+a(-1)c(-1)[b,d]+a(-1)b(-1)[c,d]\Big).
\end{multline}
All the terms in the expansion have the form $u(-1)v(-1)w$ for some 
$u,v,w\in V_1$.  The zero mode of this form is $\big(u(-1)v(-1)w\big)(2)$.
We expand $(u(-1)v(-1)w)(2)$ as an operator on $V_1$ by applying 
associativity twice.
\begin{multline*}
(u(-1)v(-1)w)(2) =\\
u(-1)w(1)v(0)+v(-1)w(1)u(0)+w(0)v(0)u(0)+w(-1)v(1)u(0).
\end{multline*}
Let $x$ be in $V_1$ and compute that $w(-1)v(1)u(0)x = 
\langle v|[u,x]\rangle w$.
Therefore 
\[\tr\Big|_{V_1} w(-1)v(1)u(0) = \langle v|[u,w]\rangle.\]
Putting this together with the expansion of $(u(-1)v(-1)w)(2)$, we get
\begin{equation}
\label{u(-1)v(-1)w}
\tr\Big|_{V_1}(u(-1)v(-1)w)(2) 
= -\langle [u,v]|w \rangle+\tr\Big|_{V_1}w(0)v(0)u(0)).
\end{equation}

Using equation~\eqref{u(-1)v(-1)w} on each term in 
expression~\eqref{3-terms_expan}, we compute the contribution
of the weight-three terms to $\tr\big|_{V_1}o(a[-1]b[-1]c[-1]d)$,
\begin{multline*}
\frac{1}{2}\bigg(-\langle [[a,b],c]|d \rangle -\langle [b,[a,c]]|d\rangle
-\langle [b,c]|[a,d]\rangle -\langle [a,[b,c]]|d\rangle -\langle [a,c]|[b,d]
\rangle\\
-\langle [a,b]|[c,d]\rangle
+\tr\Big|_{V_1}\Big(d(0)c(0)[a,b](0)+d(0)[a,c](0)b(0)+[a,d](0)c(0)b(0)\\
+d(0)[b,c](0)a(0)+[b,d](0)c(0)a(0)+[c,d](0)b(0)a(0)\Big)\bigg).
\end{multline*}
To simplify this expression we use the properties of $\langle\ |\ \rangle$
including equation~\eqref{[a,c]id} on the first six terms;  for the 
remaining terms, we expand the commutators using $[u,v](0)=u(0)v(0)-v(0)u(0)$ 
and then use $\tr AB= \tr BA$ to move $a(0)$ to the left of each term. 
We find the trace of the weight-three terms on $V_1$ is
\begin{equation}
\label{tr_3-terms}
-\langle [a,b]|[c,d]\rangle+\frac{1}{2}\tr\Big|_{V_1}\big(a(0)b(0)d(0)c(0)+
a(0)c(0)d(0)b(0)-2a(0)d(0)c(0)b(0)\big).
\end{equation}

There is only one weight-four term, $a(-1)b(-1)c(-1)d$. 
We expand its zero mode as an operator on $V_1$ using associativity
three times. 
\begin{multline*}
\big(a(-1)b(-1)c(-1)d)\big)(3)=a(-1)d(1)c(0)b(0)+b(-1)d(1)c(0)a(0)\\
                 +c(-1)d(1)b(0)a(0)+d(0)c(0)b(0)a(0)+d(-1)c(1)b(0)a(0).
\end{multline*}
Now let $x$ be in $V_1$; 
$a(-1)b(1)c(0)d(0)x=\langle b|[c,[d,x]]\rangle a$.
Thus
\begin{equation*}
\tr\Big|_{V_1}a(-1)b(1)c(0)d(0)=\langle [b,c]|[d,a]\rangle.
\end{equation*}
Using this, the expansion of $\big(a(-1)b(-1)c(-1)d\big)(3)$ and 
equation~\eqref{[a,c]id}, 
we compute
\begin{align}
\nonumber
\tr\Big|_{V_1}(a(-1)b(-1)c(-1)d)(3) &= -\langle [a,b]|[c,d]\rangle 
                                       - 2 \langle [a,d]|[b,c]\rangle \\
\label{tr_4-term}
                                    &\quad+\tr\Big|_{V_1}a(0)d(0)c(0)b(0).
\end{align}
This is the contribution of the weight-four term to 
$\tr\big|_{V_1}o(a[-1]b[-1]c[-1]d)$.   

We now add the contributions of the round-bracket weights 
zero, two, three and four terms from expressions 
\eqref{tr_0-terms}, \eqref{tr_2-terms}, \eqref{tr_3-terms} and
\eqref{tr_4-term} to find that 
\begin{equation}
\label{a[-1]b[-1]c[-1]d1}
\begin{array}{rcl}
\lefteqn{\tr\Big|_{V_1}o(a[-1]b[-1]c[-1]d)=
-\frac{n+240}{720}\Big(4\langle [a,b]|[c,d]\rangle-5\langle [a,d]|[b,c]\rangle
\Big)}\hspace{0cm}&&\\
&&+\dfrac{1}{3}\kappa([a,b],[c,d])+\dfrac{1}{6}\kappa([a,d],[b,c])+\dfrac{1}{2}
\tr\Big|_{V_1}\big(a(0)b(0)d(0)c(0)\\
&&+a(0)c(0)d(0)b(0)\big)
+\dfrac{n-48}{1152}\sum_{\pi\in S_4}\langle\pi(a)|\pi(b)\rangle
\langle \pi(c)|\pi(d)\rangle\\
&&-\dfrac{1}{48}\sum_{\pi\in S_4}\langle \pi(a)|\pi(b)\rangle \kappa(\pi(c),
\pi(d))
\end{array}
\end{equation}
(Recall that the trace of the weight-one terms on $V_1$ is zero.)
We uncover another symmetry in this trace by expanding 
$\kappa([a,b],[c,d])$ and $\kappa([a,d],[b,c])$.
\begin{align*}
\kappa([a,b],[c,d])&=\tr\Big|_{V_1}[a,b](0)[c,d](0)\\
                   &=\tr\Big|_{V_1}\big(a(0)b(0)c(0)d(0)-a(0)b(0)d(0)c(0)\\
                   &\quad-a(0)c(0)d(0)b(0)+a(0)d(0)c(0)b(0)\big)
\end{align*}
Likewise, 
\begin{align*}
\kappa([a,d],[b,c]) &= \tr\Big|_{V_1}\big(a(0)d(0)b(0)c(0)-a(0)d(0)c(0)b(0)\\
                    &\quad  -a(0)b(0)c(0)d(0)+a(0)c(0)b(0)d(0)\big).
\end{align*}
So, 
\begin{multline*}
\frac{1}{3}\kappa([a,b],[c,d])+\frac{1}{6}\kappa([a,d],[b,c])+\frac{1}{2}
\tr\Big|_{V_1}\big(a(0)b(0)d(0)c(0)+a(0)c(0)d(0)b(0)\big)\\
=\frac{1}{6}\sum_{\pi\in S_3}a(0)\pi\big(b(0)\big)\pi\big(c(0)\big)
\pi\big(d(0)\big).
\end{multline*}
Substituting this into equation~\eqref{a[-1]b[-1]c[-1]d1} completes the 
proof of the second equation in the lemma.
\end{proof}

In order to define elements of $V_{[4]}$ with zero trace on $V_0$ and 
non-zero trace on $V_1$, we need to know about the structure of $V_1$.
We will use results from Dong and Mason's paper \textit{Holomorphic 
vertex operator algebras of small central charge}, which we summarize 
here.~\cite{DM:hVOAscc} 

Dong and Mason studied the structure of strongly-rational holomorphic VOAs 
with central charges, $c$, equal to $8$, $16$ or $24$. As we have already 
completed the proof of Theorem~\ref{main_c=8,16} in the $c=8$ case, 
we need to know about $c=16$ and $c=24$.  For $c=16$, Dong and Mason show
that we have a lattice VOA, $V_L$, where $L$ is one of the two 
unimodular rank 16 lattices $E_8+E_8$ or $D_{16}^+$.  
The 
$E_8+E_8$ lattice VOA has $V_1$ isomorphic to the Lie algebra of type 
$E_8\oplus E_8$ and the $D_{16}^+$ lattice VOA has $V_1$ isomorphic
to the Lie algebra of type $D_{16}$.  For $c=24$, Dong and Mason find that
$V_1$ is either $0$, abelian of rank $24$ or semi-simple of rank less than
or equal to $24$. 

We will deal with the two $c=16$ cases and all of the $c=24$ cases except 
$V_1=0$, $V_1$ abelian or $V_1$ simple of type $A_1$, $A_2$, $D_4$, $E_6$,
$E_7$, $E_8$, $G_2$, or $F_4$ using 
Lemmas~\ref{semisimple-simple} and~\ref{simple-exceptional}.
Because we will frequently refer to the list of exceptional Lie algebras 
plus $A_1$, $A_2$ and $D_4$, we will call this list the \emph{augmented 
exceptional list} of Lie algebras.
We address $V_1$ on the augmented exceptional 
list using Lemma~\ref{exceptional} and $V_1$ abelian using Lemma~\ref{abelian}.
We do not consider the $V_1=0$ case, as our method uses elements of $V_1$.  
Dong and Mason describe the space of graded traces for
the Moonshine module, for which $V_1=0$.~\cite{DM:MM}  It is conjectured
that this is the only $c=24$, holomorphic VOA with $V_1=0$.~\cite{FLM}.

\begin{remark}
This augmented exceptional list of Lie algebras is the same list considered
by Deligne in his paper ``La S\'erie exceptionnelle de groupes de Lie.''
~\cite{De:SegL}
\end{remark}

The elements of $V_{[4]}$ used in Lemmas~\ref{semisimple-simple} 
and~\ref{simple-exceptional} are defined 
to have zero trace on $V_0$.  
We start by choosing $u$ and $v$ in $V_1$, orthogonal and of the same length 
with
respect to the Li-Zamilodichov metric and commuting with respect to the 
Lie bracket.  We then define $x(u,v)\in V_{[4]}$:
\begin{equation}
\label{x(u,v)_defn}
x(u,v):= u[-1]^3u-6u[-1]^2v[-1]v+v[-1]^3v.
\end{equation}
It is now a straight-forward application of Lemma~\ref{general_traces} to 
compute the traces of $x(u,v)$ on $V_0$ and $V_1$.
\begin{lemma}
\label{tr_x(u,v)}
Let $V$ be a strongly rational VOA, let $u$ and 
$v$ be orthogonal 
commuting elements of $V_1$ such that 
$\langle u|u\rangle =\langle v|v\rangle$, and let $x(u,v)$ be as defined
in equation~\eqref{x(u,v)_defn} above.
Then 
\begin{align*}
\tr\big|_{V_0}o(x(u,v))&=0,\\
\tr\big|_{V_1}o(x(u,v))&=\tr\big|_{V_1}\big(u(0)^4-6u(0)^2v(0)^2+v(0)^4\big).
\end{align*}
\end{lemma}  

In the case that $V_1$ is semisimple but not on the augmented exceptional 
list,
we show that there exist  $u$ and $v$ satisfying the conditions of 
Lemma~\ref{tr_x(u,v)} such that the trace of $o(x(u,v))$ on $V_1$ is 
nonzero.  
We choose linearly independent $u$ and $v$ in 
a Cartan subalgebra, $\mathfrak{h}$, of $V_1$, which is possible unless $V_1$ 
is simple of type $A_1$. Such $u$ and $v$ commute and it is easy to compute
$\tr|_{V_1}(u(0)^4-6u(0)^2v(0)^2+v(0)^4)$.  Indeed, suppose $\mathfrak{h}$ 
has root system $\Phi$ and root spaces $L_\alpha$.  The Cartan 
decomposition of $\mathfrak{g}$ is: $V_1=\mathfrak{h}
\oplus\sum_{\alpha\in\Phi}L_\alpha$.  Clearly, for any $u\in\mathfrak{h}$, 
$u(0)x=[u,x]=0$ for any $x$ in $\mathfrak{h}$, while
$u(0)x=[u,x]=\alpha(u)x$ for any $x$ in $L_\alpha$.  As each $L_\alpha$ is 
one dimensional, it follows that: For any $u$ and $v$ in the Cartan 
subalgebra of $V_1$,
\begin{equation}
\label{root_sum}
\tr\big|_{V_1}(u(0)^4-6u(0)^2v(0)^2+v(0)^4)
=\sum_{\alpha\in\Phi}\alpha(u)^4-6\alpha(u)^2\alpha(v)^2+\alpha(v)^4
\end{equation}
 
In addition to commuting, we need $u$ and $v$ to be orthogonal and of the 
same length. 
The semisimple Lie algebra $V_1$ has a decomposition into simple Lie 
algebras.
\begin{equation}
\label{V_1_decomp}
V_1 = \mathfrak{g}_1\oplus\mathfrak{g}_2\oplus 
\cdots \oplus\mathfrak{g}_m,
\end{equation}
Such a decomposition is orthogonal with respect to any invariant bilinear form,
in particular the Li-Zamilodichov metric. 

In our first case, $V_1$ semi-simple but not simple, i.e., $m>1$. 
We guarantee the orthogonality of $u$ and $v$ by choosing 
$u\in\mathfrak{g}_1$ and $v\in\mathfrak{g}_2$. 
We will actually pick $u$ and $v$ in real subspaces of $\mathfrak{h}_1$ and
$\mathfrak{h}_2$, the Cartan subalgebras of $\mathfrak{g}_1$ and 
$\mathfrak{g}_2$, respectively.  Let $\Phi_i$ be the root system of 
$\mathfrak{g}_i$.  For each root $\alpha\in\Phi_i$, fix a 
non-zero $e_\alpha\in L_\alpha$. 
There exist
$f_\alpha\in L_{-\alpha}$ and 
$h_\alpha:=[e_\alpha,f_\alpha]\in\mathfrak{h}_i$ such that 
$[h_\alpha,e_\alpha]=2e_\alpha$ and $[h_\alpha,f_\alpha]=-2f_\alpha$.   
Define 
$\mathfrak{r}_i=\Span_{\mathbb{R}}\{h_\alpha|\alpha\in\Phi_i\}$.
The  Killing form is a real-valued positive-definite 
form on 
$\mathfrak{r}_i$.  Furthermore, Dong and Mason prove that, for 
strongly rational, holomorphic VOAs 
with central charge at most 24, 
\begin{equation}
\label{kappa_to_L-Z}
\kappa(u,v)=2\langle u|v\rangle\big( \frac{n}{c}-1\big),
\end{equation}
where $n=\dim V_1$~\cite{DM:hVOAscc}.
In particular, as $c$ is an integer and for $V_1$ semisimple 
both forms are nondegenerate on $V_1$, 
$\langle\ |\ \rangle$ is a non-zero rational multiple of the 
Killing form. Hence it is also positive definite on each $\mathfrak{r}_i$. 
Thus we can choose 
$u$ in $\mathfrak{r}_1$ and $v$ in $\mathfrak{r}_2$ such that
$\langle u|u\rangle =\langle v|v\rangle\ne 0$. 
\begin{lemma}  
\label{semisimple-simple}
Let $V$ be a strongly rational holomorphic VOA with $c\le 24$,
and suppose that $V_1$ is semisimple, but not simple.  
Then there exist $u$ and $v$ in $V_1$ such that
\begin{align*}
\tr\big|_{V_0}o(x(u,v))&=0,\\
\tr\big|_{V_1}o(x(u,v))&\ne 0.
\end{align*}
\end{lemma}
\begin{proof} From the above discussion there exist $u\in\mathfrak{r}_1$ 
and $v\in\mathfrak{r}_2$ such that 
$\langle u|u\rangle =\langle v|v\rangle \ne 0$.
As $\mathfrak{r}_i\subset\mathfrak{g}_i$, $u$ and $v$ are commuting and 
orthogonal.  Thus $u$ and $v$ satisfy the conditions of 
Lemma~\ref{tr_x(u,v)}.  It follows immediately that 
$\tr\big|_{V_0}o(x(u,v))=0$, and, since $\mathfrak{r}_i\subset\mathfrak{h}_i$,
equation~\eqref{root_sum} implies that 
\[
\tr\big|_{V_1}o(x(u,v))=\sum_{\alpha\in\Phi}\alpha(u)^4
                         -6\alpha(u)^2\alpha(v)^2+\alpha(v)^4,
\]
where $\Phi$ is a root system for $V_1$.
Let $\Phi_i$ be the root system for $\mathfrak{g}_i$. Then
$\Phi = \Phi_1+\Phi_2+\cdots+\Phi_m$, and $\alpha(u)=0$, for any 
$\alpha\notin\Phi_1$, likewise, $\alpha(v)=0$ for any $\alpha\notin\Phi_2$.  
Thus
\[ 
\tr\big|_{V_1}o(x(u,v))=\sum_{\alpha\in\Phi_1}\alpha(u)^4
+\sum_{\alpha\in\Phi_2}\alpha(v)^4.
\]
For all roots $\alpha$ and $\beta$, $\alpha(h_\beta)$ is an integer. Since
$u$ and $v$ are in the real span of the $h_\beta$'s, $\alpha(u)$ and 
$\alpha(v)$ are both real numbers. Hence, 
$\tr\big|_{V_1}o(x(u,v))\ge 0$. Furthermore, because
$u\ne 0$ and $\Phi_1$ spans $\mathfrak{h}_1^*$, there exists an 
$\alpha\in\Phi_1$ such that $\alpha(u)\ne 0$. Hence 
$\tr\big|_{V_1}o(x(u,v))>0.$  
\end{proof}

We are left with the cases $V_1$ abelian and $V_1$ simple.  
If $V_1$ is abelian, then we cannot use Lemma~\ref{tr_x(u,v)} to get
$\tr\big|_{V_1}o(x(u,v))\ne 0$, but we can use this lemma to address
all the $V_1$ simple cases except those on the 
augmented exceptional list.  

Assume $V_1$ is simple, but not on the augmented exceptional list. 
We choose $u$ and $v$ in the Cartan subalgebra of $V_1$ so that $[u,v]=0$.
We also need $u$ and $v$ orthogonal and of equal length with respect to 
the form
$\langle\ |\ \rangle$.  As we mentioned previously, for $V$ strongly rational 
and holomorphic with $c\le 24$ and $V_1$ semisimple, $\langle\ |\ \rangle$ is 
a non-zero rational multiple of the Killing form on $V_1$. 
As $V_1$ is not $A_1$ type, the dimension of the Cartan subalgebra is greater 
than one, and we can choose $u$ and $v$ orthogonal and of the
same length with respect to the Killing form. They will have the same 
properties with respect to $\langle\ |\ \rangle$.  
\begin{lemma}
\label{simple-exceptional}
Let $V$ be a strongly rational VOA and suppose that $V_1$ is a simple Lie 
algebra of type $A$, $B$, $C$ or $D$ but not
type $A_1$, $A_2$ or $D_4$.  There exist $u$ and $v$ in $V_1$ such that 
\begin{align*}
\label{simple_tr_V_0}
\tr\big|_{V_0}o(x(u,v))&=0\\
\tr\big|_{V_1}o(x(u,v))&\ne 0
\end{align*}
\end{lemma}
\begin{proof} For each root type we pick explicit $u$ and $v$ in the 
Cartan subalgebra  $\mathfrak{h}$ of $V_1$, orthogonal and of equal length 
with respect to the Killing form.  From our discussion above, such a 
$u$ and $v$ satisfy the conditions of Lemma~\ref{tr_x(u,v)}, thus the trace of 
$o(x(u,v))$ on $V_0$ is 0.  

Let $(\ ,\ )$ be the nondegenerate bilinear 
form induced on $\mathfrak{h}^*$ by the Killing form.  Then for 
$u$ in $\mathfrak{h}$ and $\alpha$ in $\mathfrak{h}^*$, 
$\alpha(u)=(\alpha,u^*)$.  Using this, Lemma~\ref{tr_x(u,v)} and 
equation~\eqref{root_sum}, we 
have:  
\begin{equation}
\label{inner_prod_tr}
\tr\big|_{V_1}o(x,(u,v))=\sum_{\alpha\in\Phi}(\alpha,u^*)^4-
6(\alpha,u^*)^2(\alpha,v^*)^2+(\alpha,v^*)^4,
\end{equation}
where $\Phi$ is the root system for $V_1$.
Letting $\{e_1,e_2,\ldots,e_\ell\}$ be an orthonormal set, 
the root systems of type $A$, $B$, $C$ and $D$ can be 
realized as follows.
\[
\begin{array}{|c|c|c|}\hline
       & &\text{Root system}\\ \hline
{\rule[-6pt]{0cm}{18pt}A_\ell} &\ell\ge 1
     & \{\pm(e_i-e_j)|1\le i<j \le \ell +1\}\\ \hline
{\rule[-6pt]{0cm}{18pt}B_\ell} &\ell\ge 2&\{\pm e_i, \pm(e_i\pm e_j)|1\le i<j 
     \le \ell\}\\ \hline
{\rule[-6pt]{0cm}{18pt}C_\ell} &\ell\ge 3& \{\pm 2e_i, \pm(e_i\pm e_j)|1\le i<j
     \le \ell\}\\ \hline
{\rule[-6pt]{0cm}{18pt}D_\ell} &\ell\ge 4& \{\pm(e_i\pm e_j)|1\le i<j \le 
      \ell\}\\ \hline
\end{array}
\]
For $A_\ell$, $\ell\ge 3$, choose $u$ so that $u^*=(e_{1}-e_{2})$ 
and $v$ so that $v^*=(e_{\ell}-e_{\ell+1})$.
Since $\ell \ge 3$, $u$ and $v$ are perpendicular. 
Computing using equation~\eqref{inner_prod_tr},
\[
\tr\big|_{A_\ell(\mathbb{C})}o(x(u,v))= 8\ell +8.
\] 
So for $\ell\ge 3$, $\tr\big|_{A_\ell(\mathbb{C})}o(x(u,v))\ne 0$.
For $B_\ell$, $C_\ell$ and $D_\ell$, 
we choose $u$ so that $u^*=e_1$ and $v$ so that $v^*=e_\ell$.  
From equation~\eqref{inner_prod_tr},
\begin{align*}
\tr\big|_{B_\ell(\mathbb{C})}o(x(u,v))&=8\ell-28,\\
\tr\big|_{C_\ell(\mathbb{C})}o(x(u,v))&=8\ell+32,\\
\tr\big|_{D_\ell(\mathbb{C})}o(x(u,v))&=8\ell-32.
\end{align*} 
Therefore $\tr\big|_{B_\ell(\mathbb{C})}\upsilon(u,v)\ne 0$ 
for any integer $\ell$, 
$\tr\big|_{C_\ell(\mathbb{C})}\upsilon(u,v)\ne 0$ for any positive $\ell$ and 
$\tr\big|_{D_\ell(\mathbb{C})}\upsilon(u,v)\ne 0$ for $\ell \ne 4$.
\end{proof}

\begin{remark}
One can actually show that for the $V_1$ on the augmented exceptional list, 
it is not 
possible to choose $u$ and $v$ as in Lemma~\ref{simple-exceptional} so that 
the trace of $x(u,v)$ on $V_1$ is non-zero.  For a proof see~\cite{KLH:thesis}.
In the type $A$, $D$ and $E$ cases, this reflects the fact that $A_1$, $A_2$,
$D_4$, $E_6$, $E_7$ and $E_8$ are the strongly perfect root lattices, i.e.,
those root lattices which are spherical 4-designs.~\cite{Mar:ples}
\end{remark}

Both $c=16$ cases; $V_1$ of type $E_8+E_8$ and $V_1$ of type $D_{16}$ have 
been addressed, by Lemmas~\ref{semisimple-simple} and~\ref{simple-exceptional}
respectively.  Thus we are left with the $V_1$ on the augmented exceptional 
list with $c=24$

Recently, in their paper \emph{Integrability of $C_2$-cofinite vertex operator 
algebras}, Dong and Mason prove that the levels of the affine Lie algebras 
created by the simple subalgebras of $V_1$ must be positive 
integers.~\cite{DM:IVOA}
We do not define level here, but once one understands it, it is a relatively
simple computation to show that for $c=24$ VOAs with $V_1$ on the augmented
exceptional list, only $D_4$-type $V_1$ has an integral level.  
To do this one uses the 
values in table~\ref{Lie_alg_calcs}.  Thus most of the cases we will address 
in Lemma~\ref{exceptional} are actually moot. In fact, as a VOA with $V_1$ 
of $D_4$ type does not appear on Schellekens list of holomorphic $c=24$ 
VOAs~\cite{Schell}, 
the $D_4$ case is likely moot as well.  However, as we do not 
want to rely on Schellekens list, which is not mathematically rigorous, we
proceed with the $D_4$ case.  And, since it merely adds some lines to 
table~\ref{Lie_alg_calcs}, which we want to include to aid the reader in his 
or her computation of the levels, we proceed with the other 
augmented-exceptional cases as well.

To  address the $c=24$ VOAs with $V_1$ on the augmented exceptional list, 
we need a different
approach to that taken in Lemmas~\ref{semisimple-simple} 
and~\ref{simple-exceptional}. Recall 
equation~\eqref{kappa_to_L-Z}, relating the Killing form to the 
Li-Zamilodichov metric, 
for strongly rational holomorphic VOAs with central charge less than 24 and 
$V_1$ semisimple. For $c=24$, it yields
\begin{equation}
\label{angle_to_killing}
\langle u|v\rangle =\frac{12}{n-24}\kappa(u,v),
\end{equation}
where $n$ is the dimension of $V_1$.  (Note that this implies that $n\ne 24$ 
for $V_1$ semi-simple.  In fact, $n=24$ if and only if $V_1$ is 
abelian~\cite{DM:hVOAscc}.)  
Using equation~\eqref{angle_to_killing} to 
rewrite Lemma~\ref{tr_x(u,v)} in terms of the Killing form yields the following
lemma.
\begin{lemma}
\label{general_traces_kappa} 
Let $V$ be a strongly rational holomorphic $c=24$ VOA with $V_1$ 
semisimple, and let $a$, $b$, $c$ and $d$ be in $V_1$.  Then
\begin{align*}
\tr\big|_{V_0} o(a[-1]b[-1]c[-1]d) &=
    -\frac{1}{60(n-24)}\big(4\kappa([a,b],[c,d])+5\kappa([a,d],[b,c])\big)\\
& \quad +\frac{1}{8(n-24)^2}\sum_{\pi\in S_4}
\kappa(\pi(a),\pi(b))\kappa(\pi(c),\pi(d)),\\
\tr\big|_{V_1} o(a[-1]b[-1]c[-1]d) &=
 -\frac{n+240}{60(n-24)}\big(4\kappa([a,b],[c,d])+5\kappa([a,d],[b,c])\big)\\
& \quad -\frac{n}{8(n-24)^2}\sum_{\pi\in S_4}
\kappa(\pi(a),\pi(b))\kappa(\pi(c),\pi(d))\\
& \quad +\frac{1}{6}\sum_{\pi\in S_3}
\tr\Big|_{V_1}a(0)\pi\big(b(0)\big)\pi\big(c(0)\big)\pi\big(d(0)\big).
\end{align*}
\end{lemma}

For each root $\alpha$ of $V_1$ there exist $e\in L_\alpha$, 
$f\in L_{-\alpha}$ and $h\in \mathfrak{h}$ satisfying the standard 
$\mathfrak{sl}_2$ commutation relations:
\begin{align*}
 [e,f]&=h& [h,e]&=2e& [h,f]&=-2f.
\end{align*}
From Lemma~\ref{general_traces_kappa}, we compute
\begin{align}
\label{hhhh_trace_on_V_0}
\tr\Big|_{V_0}o(h[-1]^3h)&=\frac{3\kappa(h,h)^2}{(n-24)^2},\\
\label{hhhh_trace_on_V_1}
\tr\Big|_{V_1}o(h[-1]^3h)&=-\frac{3n\kappa(h,h)^2}{(n-24)^2}
       +\tr\big|_{V_1}h(0)^4,
\end{align}
and
\begin{align}
\label{efef_trace_on_V_0}
\tr\Big|_{V_0}o(e[-1]f[-1]e[-1]f)&=\frac{\kappa(h,h)}{60(n-24)}
                                  +\frac{\kappa(h,h)^2}{2(n-24)^2},\\
\label{efef_trace_on_V_1}
\tr\Big|_{V_1}o(e[-1]f[-1]e[-1]f)&=
\frac{\kappa(h,h)(n+240)}{60(n-24)}-\frac{n\kappa(h,h)^2}{2(n-24)^2}
+\frac{1}{6}\tr\big|_{V_1}h(0)^4.
\end{align}
For the traces of $o(e[-1]f[-1]e[-1]f)$ we use $2\kappa(e,f)=\kappa(h,h)$,
which follows easily from the invariance of the Killing form, the 
commutation relations and 
\begin{equation}
\label{efef_to_hhhh}
\tr\big|_{V_1}\big(2e(0)f(0)e(0)f(0)+4e(0)^2f(0)^2\big)=
\tr\big|_{V_1}h(0)^4.
\end{equation}
To prove equation~\eqref{efef_to_hhhh}, we decompose $V_1$ as a module 
for $\mathfrak{s}_\alpha$, the subalgebra generated by $e$, $f$ 
and $h$, 
$V_1=M_1\oplus M_2 \oplus \cdots \oplus M_m$.  Each $M_i$ is an irreducible 
$\mathfrak{sl}_2$-module with dimension at most 4.  For each dimension, such
a module is unique up to isomorphism and  
completely known. So we can just check that 
$\tr\big|_{M_i}\big(2e(0)f(0)e(0)f(0)+4e(0)^2f(0)^2\big)=
\tr\big|_{M_i}h(0)^4$ holds for each of them.

Now define $C$ to be
$\frac{\tr|_{V_0}o(e[-1]f[-1]e[-1]f)}{\tr|_{V_0}o(h[-1]^3h)}$.
From equations~\eqref{efef_trace_on_V_0} and~\eqref{hhhh_trace_on_V_0},
\[
C=\frac{n-24}{180\kappa(h,h)}+\frac{1}{6}.
\]
Define $y(\alpha)$ to be
\begin{equation}
\label{defn_y(alpha)}
y(\alpha)=e[-1]f[-1]e[-1]f-Ch[-1]^3h.
\end{equation}
The constant $C$ is defined so that
\[\tr\big|_{V_0}o(y(\alpha))=0.\]
It is straight forward to compute the trace of $o(y(\alpha))$ on $V_1$ using
equations~\eqref{hhhh_trace_on_V_1} and~\eqref{efef_trace_on_V_1};
\begin{equation}
\label{y_alpha_on_V_1}
\tr\big|_{V_1}o(y(\alpha))=\frac{\kappa(h,h)(n+120)}
{30(n-24)}-\frac{n-24}{180\kappa(h,h)}\tr\big|_{V_1}h(0)^4.
\end{equation}

We wish to show that, for $c=24$ VOAs with $V_1$ on the augmented exceptional
list, there 
exists a root $\alpha$ such that $\tr\big|_{V_1}o(y(\alpha))\ne 0$.  
To simplify the computation, we 
choose $\alpha$ to be a long root, in which case
\begin{equation}
\label{hhhh_to_kappa(h,h)}
\tr\big|_{V_1}h(0)^4=\kappa(h,h)+24.
\end{equation}
To show this we again decompose $V_1$ into irreducible modules for 
$\mathfrak{s}_\alpha$
\[V_1=\mathfrak{s}_\alpha \oplus M_1\oplus M_2 \oplus \cdots \oplus M_k.\]
Recall that $\kappa(h,h)=\tr\big|_{V_1}h(0)^2$.  This trace, like the
trace of $h(0)^4$, splits over the decomposition.
Since $\alpha$ is a long root,  $\dim M_i\le 2$ for all $i$, and   
we compute that
\begin{align*}
\tr\big|_{M_i}h(0)^2&=\tr\big|_{M_i}h(0)^4,\\
\tr\big|_{\mathfrak{s}_\alpha}h(0)^2+24&=\tr\big|_{\mathfrak{s}_\alpha}h(0)^4.
\end{align*}  
Equation~\eqref{hhhh_to_kappa(h,h)}
follows.

Inserting equation~\eqref{hhhh_to_kappa(h,h)} 
into equation~\eqref{y_alpha_on_V_1} and 
factoring, we see that, 
for a long root $\alpha$, $\tr\big|_{V_1}o(y(\alpha))=0$ if and only if 
\[
\big(6\kappa(h,h)-n+24\big)\big((n+120)\kappa(h,h)+24(n-24)\big)=0.
\]
Since $n$ must be positive and $\kappa(h,h)\ge 8$, the second factor 
cannot be zero. Thus the trace of $o(y(\alpha))$ on $V_1$ is zero if and only 
if the first factor is zero.
Table~\ref{Lie_alg_calcs} gives the values of $n$ and $\kappa(h,h)$ for 
the on the augmented list of exceptional Lie algebras.  
Clearly no pair of values satisfies $n=6\kappa(h,h)+24$.  
\begin{table}[ht]
\label{Lie_alg_calcs}
\begin{center}
\begin{tabular}{|c|c|c|}\hline
$V_1$ & $\dim V_1$ & $\kappa(h,h)$\\ \hline
$A_1$ & $3$        & $8$          \\ \hline
$A_2$ & $8$        & $12$         \\ \hline
$D_4$ & $28$       & $24$         \\ \hline
$E_6$ & $78$       & $48$         \\ \hline
$E_7$ & $133$      & $72$         \\ \hline
$E_8$ & $248$      & $120$        \\ \hline
$F_4$ & $52$       & $36$         \\ \hline
$G_2$ & $14$       & $16$         \\ \hline
\end{tabular}
\end{center}
\end{table}
This gives us the following lemma.
\begin{lemma}
\label{exceptional}
Let $V$ be a strongly rational homomorphic $c=24$ vertex operator algebra
 with $V_1$ simple
of type $A_1$, $A_2$, $D_4$ or exceptional.  Then for a long root $\alpha$ of
$V_1$,
\begin{align*} 
\tr\big|_{V_0}o(y(\alpha))&=0,\\
\tr\big|_{V_1}o(y(\alpha))&\ne 0.
\end{align*}
See equation~\eqref{defn_y(alpha)} for the definition of $y(\alpha)$.
\end{lemma}

We turn to the $V_1$ abelian case.  Dong and Mason show that in this case
$V$ is the lattice VOA generated by the Leech lattice.~\cite{DM:hVOAscc}
As the Killing form on an abelian 
Lie algebra is zero, while the Li-Zamilodichov metric is non-degenerate 
on $V_1$, equation~\eqref{kappa_to_L-Z} implies that, for strongly rational
 holomorphic 
$c=24$ VOAs with $V_1$ abelian, $\dim V_1= n =24$.  Let $a$ be an element of
$V_1$, and rewrite $a[-1]^3a$ as $a[-1]^4\boldsymbol{1}$.  From
Lemma~\ref{general_traces}
\begin{equation}
\label{traces_a[-1]^4}
\begin{split}
\tr\big|_{V_0}o(a[-1]^4\boldsymbol{1})&=\frac{1}{48}\langle a|a\rangle^2, \\
\tr\big|_{V_1}o(a[-1]^4\boldsymbol{1})&=-\frac{1}{2}\langle a|a\rangle^2. 
\end{split}
\end{equation}

Next we compute the traces of $a[-2]^2\boldsymbol{1}$.
For this we use the following equation of 
Zhu:
\[ 
Z(a[-2]b,\tau)=-\sum_{k=2}^\infty (2k-1)E_{2k}(\tau)Z(a[2k-2]b,\tau)
\]
\cite[Prop. 4.3.6]{Zhu:mi}; see also~\cite[Prop 2.2.1]{DMN}.
Let $b=a[-2]\boldsymbol{1}$.  As $a[2k-2]a[-2]\boldsymbol{1}=0$ for 
$k\ge 3$, the above equation becomes
\[
Z(a[-2]^2\boldsymbol{1},\tau)=-3E_4(\tau)Z(a[2]a[-2]\boldsymbol{1},\tau).
\]
Furthermore, from the affine Lie algebra commutation relations
$a[2]a[-2]\boldsymbol{1}= 2\langle a|a \rangle\boldsymbol{1}$, so 
\begin{align*}
Z(a[-2]^2\boldsymbol{1},\tau)
&=-6 \langle a|a \rangle E_4(\tau)Z(\boldsymbol{1},\tau),\\
&=-6 \langle a|a \rangle (\frac{1}{720}+\frac{1}{3}q +\cdots)
                         (q^{-1}+24+\cdots),\\
&= -\frac{\langle a|a \rangle}{120}q^{-1}-\frac{11\langle a|a\rangle}{5}+
\cdots.
\end{align*}
Thus 
\begin{equation}
\label{traces_a[-2]^2}
\begin{split}
\tr\big|_{V_0}o(a[-2]^2\boldsymbol{1})&=-\frac{\langle a|a \rangle}{120},\\
\tr\big|_{V_1}o(a[-2]^2\boldsymbol{1})&=-\frac{11\langle a|a\rangle}{5}.
\end{split}
\end{equation}
\begin{lemma}
\label{abelian}
Let $V$ be a strongly holomorphic $c=24$ VOA with $V_1$ abelian. 
Then there exists an element $a$ in $V_1$ such that  
\begin{align*}
\tr\big|_{V_0}o((2a[-1]^4+5\langle a|a\rangle a[-2]^2)\boldsymbol{1})&=0,\\
\tr\big|_{V_1}o((2a[-1]^4+5\langle a|a\rangle a[-2]^2)\boldsymbol{1})&\ne 0.
\end{align*}
\end{lemma}
\begin{proof}  Since the Li-Zamilodichov metric is nondegenerate on $V_1$, 
there exists an $a\in V_1$ such that $\langle a|a\rangle\ne 0$. 
From equations~\eqref{traces_a[-1]^4} and~\eqref{traces_a[-2]^2}, 
we compute that
\begin{align*}
\tr\big|_{V_0}o((2a[-1]^4+5\langle a|a \rangle a[-2]^2)\boldsymbol{1})&=0,\\
\tr\big|_{V_1}o((2a[-1]^4+5\langle a|a \rangle a[-2]^2)\boldsymbol{1})&=
-12\langle a|a\rangle^2.
\end{align*}  
\end{proof}

This completes the last case in the proof of the main theorems.

\section{Conclusion}

Having seen that for strongly rational holomorphic VOAs of small central 
charge there are no nontrivial obstructions to attaining the expected modular 
forms as graded traces, we can explore the question of which modular forms can
be attained as the graded traces of highest-weight vectors.  Dong, Mason and 
Nagatomo conjecture that there are no nontrivial obstructions in that case 
either.~\cite{DMN}  They prove their conjecture in the $c=8$ case and produce
an interesting class of highest-weight vectors for lattice VOAs using 
spherical harmonics and compute 
their graded traces.~\cite{DMN}  For the moonshine module, Dong and Mason 
find another class of highest-weight vectors with weight divisible by four and 
compute their graded traces.~\cite{DM:MM}.  For the same VOA, the author finds 
a related family of highest-weight vectors and their graded 
traces.~\cite{KLH:hwmm}  Together these papers show that for the moonshine 
module all cusp forms with weights 12, 16, 18, 20, 22 and 26 can be attained
as graded traces.  (In the moonshine module cases, $V_1=0$ so the forms must 
be cusp forms.)

Dong, Mason and Nagatomo's use of spherical harmonics to produce highest-weight
vectors for lattice VOAs can be generalized to produce highest-weight vectors
for all VOAs of strong CFT-type.  In fact, the vector $x(u,v)$ used in 
Lemma~\ref{tr_x(u,v)} is an example of such a highest-weight 
vector.~\cite{KLH:thesis}
As in this paper, one can compute their traces on $V_0$ and $V_1$ relative to 
the
root system of $V_1$, but this only determines the graded traces when the 
weight is low enough for the relevant space of modular forms to be one 
dimensional.  For example, in the $c=24$, $V_1\ne 0$ case, 
the methods of this paper should be able to 
produce highest-weight vectors with specific graded traces for weights 
6, 8 and 10, but not for higher weights. To compute the graded traces of 
given highest-weight vectors for general VOAs is difficult as one needs 
knowledge of the structure of $V_n$ for all $n$.  
Completion of the classification of holomorphic VOAs of small central charge 
would make this project feasible in those cases.~\cite{Schell}, 
\cite{DM:hVOAscc}, \cite{DM:IVOA}.

\section{Acknowledgments}

The author would like to thank her dissertation adviser, Professor Geoffrey 
Mason, and the other members of her dissertation committee, Professors 
Chongying Dong and Michael Tuite.

\bibliography{kbiblio}
\bibliographystyle{alpha}

\end{document}